# Parking 3-sphere swimmer
## I. Energy minimizing strokes

François Alouges

CMAP, École Polytechnique,
route de Saclay,
91128 Palaiseau Cedex,
France

Giovanni Di Fratta

CMAP, École Polytechnique,
route de Saclay,
91128 Palaiseau Cedex,
France

**Abstract.** The paper is about the parking 3-sphere swimmer ($\mathsf{sPr}_3$). This is a low-Reynolds number model swimmer composed of three balls of equal radii. The three balls can move along three horizontal axes (supported in the same plane) that mutually meet at the center of $\mathsf{sPr}_3$ with angles of 120°. The governing dynamical system is introduced and the implications of its geometric symmetries revealed. It is then shown that, in the first order range of small strokes, optimal periodic strokes are ellipses embedded in $3d$ space, i.e. closed curves of the form $t \in [0, 2\pi] \mapsto (\cos t)u + (\sin t)v$ for suitable orthogonal vectors $u$ and $v$ of $\mathbb{R}^3$. A simple analytic expression for the vectors $u$ and $v$ is derived. The results of the paper are used in a second article where the real physical dynamics of $\mathsf{sPr}_3$ is analyzed in the asymptotic range of very long arms.

## 1. Introduction

The problem of swimming at low Reynolds number is a question of considerable biological and biomedical relevance which deserves also great interest from the point of view of fundamental science. Starting from the pioneering works of Taylor [Tay51], Lighthill [Lig52] and Purcell [Pur77], the problem has received a lot of attention in recent years[1]. This problem is both surprising and attractive due to the fact that one focuses on swimmers of very small size: Indeed as pointed out by Taylor [Tay51] and Purcell [Pur77] the physics of swimming at length scales of a few micrometers is very different from our common macro-scale experience. The most striking difference is reflected by the low value of the Reynolds number. The Reynolds number $Re = \rho u L / \mu$ gives an estimate for the relative importance of inertial to viscous forces for an object of typical length scale $L$ moving at speed $u$ through a Newtonian fluid of density $\rho$ and *dynamic viscosity* $\mu$. In applications the velocity $u$ can rarely exceed a few body lengths per second and therefore, if one focuses on the swimming processes taking place in a given fluid, the Reynolds number is entirely controlled by $L$. At small $L$, inertial forces are negligible and, in order to move, micro-swimmers can only exploit the viscous resistance of the surrounding fluid.

This implies that microorganisms, such as bacteria, are required to take swimming strategies completely different from those employed by larger organisms, such as fish. In particular, the observation that, in a flow regime obeying Stokes equations, a scallop cannot advance through the reciprocal motion of its valves is encapsulated by Purcell's *scallop theorem* [Pur77]. The mathematical explanation for this is in the symmetry of the Stokes equations under time reversal: whatever forward motion will be produced by closing the valves, it will be exactly cancelled by a backward motion upon reopening them. This leads to the investigation of the simplest mechanisms capable of self-propulsion at small spatial scales. By this, we mean the ability to advance by performing a cyclic shape change – a *stroke* – in the absence of external forces. Several proposals have been put forward and analysed (see, e.g. [Tay51], [Pur77], [NG04], [AGK04], [BKS03] and [LM09]).

---

[1]. See the encyclopedia article [DAL12] for an elementary introduction to the subject and the review paper [LP09] for a comprehensive list of references.

The basic problem of swimming can be stated as follows: given a periodic history log of shape changes of a swimmer (a *stroke*), predict the corresponding history of positions and orientations in space. A related question concerns its *controllability* i.e., given an arbitrary initial position and orientation, its possibility to achieve any prescribed position and orientation in space by the means of a suitable sequence of strokes. Indeed, the unusualness of low $Re$ swimming is in that, since inertia forces are negligible, reciprocal shape changes do not contribute to a net motion, and the question of controllability is far from being trivial when the swimmer is constrained, like in the case of a scallop, to few degrees of freedom to modify its shape. In that respect, the *scallop theorem* [Pur77] is precisely a result of non-controllability for swimmers having only one degree of freedom at their disposal.

Once controllability is known, it is then natural to investigate how to reach this configuration change at minimal energetic cost. This is a question of *optimal control* from the mathematical point of view, and a question of *natural selection* from the point of view of fundamental science, which has a fundamental role in the engineering design of micro-swimmers. Indeed, as pointed out in [AGK04]: «*Microbots must swim much faster than bacteria if they are to interface with the macroscopic world. A micron-size robot swimming 100 times as fast as a bacterium, at the modest speed of* 1mm *per second, has Reynolds number* $Re = \mathcal{O}(10^{-3})$ *and, since power scales like* $u^2$, *consumes* $10^4$ *more power than a bacterium*». Micro-swimmers must therefore attempt to swim as effectively as possible, and the problem we address is to look for optimal swimming styles.

In this paper, we focus on a very interesting micro-swimmer introduced in [LM09], hereinafter referred to as the *parking 3-sphere swimmer* (sPr$_3$), which can be thought as the complementary version of Purcell's three sphere rotator introduced in [DBS05]. Full controllability of sPr$_3$, as well as of a wide class of other model swimmers, has been rigorously proved in [ADH+13], where also a method to *numerically* address the optimal control problem has been presented. Analytical investigations of the optimal control problem have only been faced for micro-swimmers which, although resembling a $3d$ object, are constrained to live in a $1d$-like world because capable to perform a net displacement of their center of mass only along a given ($1d$) direction [ADL08, ADL09]. On the other hand, as in the case of sPr$_3$, when the net displacement of the swimmer can take place in a plane, a greater number of control variables is present and the analysis is more involved.

The aim of this paper is to *analytically* address the optimal control problem for sPr$_3$ in the range of *small* strokes. In that respect the main result of the paper is Theorem 12 which reveals the complete structure of the optimal stroke – in terms of mechanical energy – which produces a given displacement, both in translation and rotation. They turn out to be planar ellipses.

The rest of the paper is organized as follows: in section 2 we give both a geometric and a kinematic description of parking 3-sphere swimmer (sPr$_3$); we then introduce the control system object of this paper. In section 3 we investigate the geometric structure of the control system by exploiting the symmetries it has to satisfy due to the governing Stokes equations. In section 4 we reveal the structure of control system in the range of small strokes. Finally, section 5 is devoted to the characterization of the energy minimizing strokes.

## 2. The swimming problem as a control problem

We focus our attention on the swimmer sPr$_3$ proposed in [LM09] and [ADH+13]. The swimmer is composed of three non intersecting balls $(B_i)_{i \in \mathbb{N}_3}$ of $\mathbb{R}^3$ centered at $b_i \in \mathbb{R}^3$ and of equal radii $a > 0$. The three balls can move along three horizontal axes (i.e. all contained in the same horizontal plane) that mutually meet at the point $c \in \mathbb{R}^3$ with fixed angles of $2\pi/3$ one to another (see Figure 1). This reflects a situation where the balls are linked together by thin jacks that are able to elongate. The viscous resistance associated with these jacks is, however, neglected, and the fluid is thus assumed to fill the whole space outside the balls, i.e the open set $\mathbb{R}^3 \setminus \cup_{i=1}^3 \bar{B}_i$.

The balls do not rotate around their axes so that the shape of the swimmer is characterized by the three lengths $\xi_1, \xi_2, \xi_3$ of its arms, measured from $c$ to the center $b_i$ of each ball. However, the swimmer can freely rotate around $c$ in the horizontal plane containing the jacks. Eventually, due to the symmetries of the system, *the swimmer stays in the horizontal plane*.

Thus, the geometrical configurations assumed by the swimmer can be described by two set of variables:

- The vector of **shape variables** $\xi := (\xi_1, \xi_2, \xi_3) \in \mathcal{M} := \left(2a/\sqrt{3}, \infty\right)^3 \subseteq \mathbb{R}^3_+$ from which relative distances $(b_{ij})_{i,j \in \mathbb{N}_3}$ between the balls are obtained, the lower bound in $\mathcal{M}$ being chosen in order to avoid overlaps of the balls.
- The vector of **positions variables**, denoted by $p = (c, \theta) \in \mathbb{R}^2 \times \mathbb{R}$ which describe the global position and orientation in space of the swimmer.

More precisely, we consider the **reference** equilateral triangle (convexly) spanned by the unit vectors $z_1, z_2, z_3 \in \mathbb{R}^2$, with $z_1 := (1,0)$, $z_2 := R^\mathsf{T}(2\pi/3)z_1$, $z_3 := R(2\pi/3)z_1$ where $R(\theta)$ stands for the planar rotation of angle $\theta$ given by the matrix:

$$R(\theta) := \begin{pmatrix} \cos\theta & -\sin\theta \\ \sin\theta & \cos\theta \end{pmatrix}. \tag{1}$$

Position and orientation in the plane are described by the coordinates of the center $c \in \mathbb{R}^2$ and the angle $\theta$ that one arm, say arm number 1 (here and hence after denoted by $\|_\mathbf{1}$), makes with the fixed direction $z_1$ (cf. **Figure** 1). Therefore we place the center of the ball $B_i$ at $b_i := c + \xi_i R(\theta) z_i$.

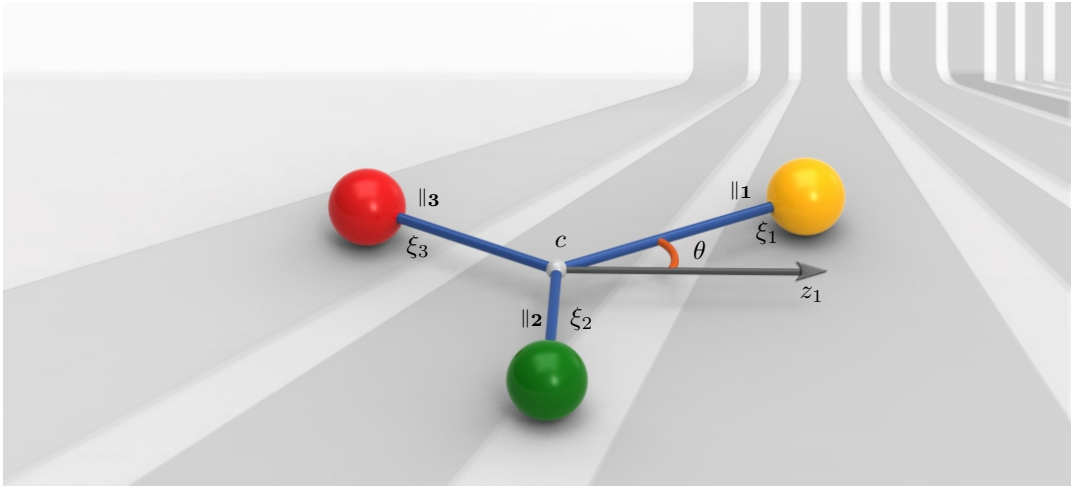

**Figure 1.** The swimmer sPr$_3$ is composed of three spheres of equal radii. The three spheres can move along three horizontal axes that mutually meet at $c$ with angle $2\pi/3$. The spheres do not rotate around their axes so that the shape of the swimmer is characterized by the three lengths $\xi_1, \xi_2, \xi_3$ of its arms, measured from the origin to the center of each ball. However, the swimmer may freely rotate around $c$ in the horizontal plane.

The swimmer is fully described by the parameters $(\xi, p) \in \mathcal{M} \times \mathbb{R}^3$. Indeed, once denoted by $B_a$ the ball of $\mathbb{R}^3$ centered at the origin and of radius $a$, for every $r \in \partial B_a$, the position of the current point on the $i$-th sphere of the swimmer in the state $(\xi, p)$ is given, for every $(\xi, p, r) \in \mathcal{M} \times \mathbb{R}^3 \times \partial B_a$, by the function

$$r_i(\xi, p, r) := c + R(\theta)(\xi_i z_i + r). \tag{2}$$

Note that the functions $(r_i)_{i \in \mathbb{N}_3}$ are analytic in $\mathcal{M} \times \mathbb{R}^3$, and we use them to compute the instantaneous velocity on the $i$-th sphere $B_i$, which for every $(\xi, p, r) \in \mathcal{M} \times \mathbb{R}^3 \times \partial B_a$ and every $i \in \mathbb{N}_3$ reads as

$$u_i(\xi, p, r) = \dot{c} + \dot{\xi}_i R(\theta) z_i + \dot{\theta} R(\theta)(\xi_i z_i^\perp + r^\perp), \tag{3}$$

with $z_i^\perp := R(\pi/2) z_i$ and $r^\perp := R(\pi/2) r$.

It has been proved in [ADH+13] that it is possible to control the state of the system sPr$_3$ (i.e. both the shape $\xi$ and the position $p$) using as controls only the rate of shape changes $\dot\xi$. To achieve this, one has to understand the way $p$ varies when one changes $\dot\xi$. This is done by assuming *self-propulsion* and that swimmer's inertia is negligible, which imply that the total viscous force and torque exerted by the surrounding fluid on the swimmer must vanish. More precisely, and we refer to [ADH+13] for the details, the control system can be written under the control form

$$\dot p = F(p \cdot e_3, \xi)\, \dot\xi \quad \text{with} \quad F(\theta, \xi) := [f_1(\theta, \xi) \,|\, f_2(\theta, \xi) \,|\, f_3(\theta, \xi)]^\mathsf{T}, \tag{4}$$

where we have used the standard notation $e_3 := (0, 0, 1) \in \mathbb{R}^3$ and denoted by $\theta := p \cdot e_3$, the angle that the arm $\|_1$ make with the fixed direction $z_1 := (0, 1)$.

Let us note that the control system $F$ does not depend on $c$ because of translational invariance of Stokes problem. On the other hand, translational invariance is just one of the symmetries which sPr$_3$ is subject to. Aim of the next section is to reveal the structure of the control system $F$ as a consequence of the symmetries it must satisfy being governed by Stokes equations.

## 3. Symmetries

For any initial condition $p_0 := (x_0, y_0, \theta_0)$ and for any control curve $\xi \colon I \subseteq \mathbb{R} \to \mathcal{M}$, with $I$ neighbourhood of zero, we denote by $\gamma(p_0, \xi) \colon I \to \mathbb{R}^3$ the solution associated to the dynamical system

$$\dot p = F(p \cdot e_3, \xi) \dot\xi \quad, \quad p(0) := p_0, \tag{5}$$

so that $\gamma(p_0, \xi)$ is such that $\forall t \in I$

$$\dot\gamma(p_0, \xi)(t) = F(\gamma(p_0, \xi)(t) \cdot e_3, \xi(t))\dot\xi(t) \quad, \quad \gamma(p_0, \xi)(0) = p_0. \tag{6}$$

### 3.1. Rotational invariance

Let us denote by $R(\theta)$ the rotation matrix that rotates by an angle $\theta$ about the $e_3 := (0, 0, 1)$ axis:

$$R(\theta) := \begin{pmatrix} \cos\theta & -\sin\theta & 0 \\ \sin\theta & \cos\theta & 0 \\ 0 & 0 & 1 \end{pmatrix}. \tag{7}$$

Rotational invariance of Stokes equations entails that the solution of the dynamical system (5) must be invariant with respect to planar rotations, i.e. that for any $\theta \in \mathbb{R}$

$$\gamma(p_0 + \theta e_3, \xi)(t) = R(\theta)(\gamma(p_0, \xi)(t) - P p_0) + \theta e_3 + P p_0, \tag{8}$$

where we have denoted by $P \colon (x, y, \theta) \mapsto (x, y, 0)$ the projection of the generic point $p$ of the state space, onto the plane generated by the vectors $e_1$ and $e_2$. Once introduced the operator $A(\theta) := [\mathrm{Id} - R(\theta)]P$, the condition of rotational invariance can be stated in the form:

CONDITION 1. (ROTATIONAL INVARIANCE) *If $\gamma(p_0, \xi)$ is a solution of the control system (5) then so is $\gamma(p_0 + \theta e_3, \xi)$ and*

$$\gamma(p_0 + \theta e_3, \xi)(t) = R(\theta)\gamma(p_0, \xi)(t) + \theta e_3 + A(\theta)p_0 \quad \forall t \in I. \tag{9}$$

By noting that $R^\mathsf{T}(\theta)A(\theta) = [R^\mathsf{T}(\theta) - \mathrm{Id}]P = -A(-\theta)$ we can write (9) in the following form which will be useful later

$$\gamma(p_0, \xi)(t) = R^\mathsf{T}(\theta)\gamma(p_0 + \theta e_3, \xi)(t) - \theta e_3 + A(-\theta)p_0 \quad \forall t \in I. \tag{10}$$

**Remark 2.** Here and in the sequel, we have preferred to state the symmetry relations satisfied by sPr₃ as hypotheses on the solution $\gamma$. The rationale behind this is in that the results work for any control system possessing the same symmetries of sPr₃, i.e. regardless of whether these hypotheses are guaranteed, as in the case of sPr₃, from the invariance of Stokes equations under a certain group of transformations.

PROPOSITION 3. *Let us denote by $\xi_0 := \xi(0) \in \mathcal{M}$ the initial state of control parameters and by $T_\xi \mathcal{M}$ the tangent space of $\mathcal{M}$ at $\xi$. If the control system (5) is invariant under rotations and for every $\xi \in \mathcal{M}$ one has $T_\xi \mathcal{M} \simeq \mathbb{R}^3$, then*

$$F(\theta, \xi) = R(\theta) F(\xi) \quad \text{with} \quad F(\xi) := F(0, \xi), \tag{11}$$

*for every $(\theta, \xi) \in \mathbb{R} \times \mathcal{M}$.*

**Proof.** Let $\gamma(p_0 + \theta e_3, \xi)$ be a solution of the control system (5). For every $\theta \in \mathbb{R}$ we have

$$\dot{\gamma}(p_0 + \theta e_3, \xi) \stackrel{(5)}{=} F(\gamma(p_0 + \theta e_3, \xi) \cdot e_3, \xi) \dot{\xi}. \tag{12}$$

On the other hand, by first using (9) and then (5), we get

$$\begin{aligned}\dot{\gamma}(p_0 + \theta e_3, \xi) &\stackrel{(9)}{=} R(\theta) \dot{\gamma}(p_0, \xi) \\ &\stackrel{(5)}{=} R(\theta) F(\gamma(p_0, \xi) \cdot e_3, \xi) \dot{\xi}.\end{aligned} \tag{13}$$

Therefore, $F(\gamma(p_0 + \theta e_3, \xi) \cdot e_3, \xi) \dot{\xi} = R(\theta) F(\gamma(x_0, \xi) \cdot e_3, \xi) \dot{\xi}$ for every $\theta \in \mathbb{R}$. Since $T_{\xi_0} \mathcal{M} \simeq \mathbb{R}^3$, $\dot{\xi}$ can be arbitrarily chosen and therefore $F(p_0 \cdot e_3 + \theta, \xi_0) = R(\theta) F(p_0 \cdot e_3, \xi_0)$ for every $\theta \in \mathbb{R}$. Eventually, referring to the initial (angular) position $p_0 \cdot e_3 = 0$ we get (11). □

### 3.2. Interchanging two arms

Let us consider the control system (5) which, due to the rotational invariance, can be written in the form (cf. (11)):

$$\dot{p} = R(p \cdot e_3) F(\xi) \dot{\xi}. \tag{14}$$

We want to investigate the effect that a swap of the arms has on the generic solution $\gamma(p_0, \xi)$. This amounts to understand how the solution behaves with respect to the action of the symmetric group $S_3$, as far as we identify its elements with the arms of sPr₃. Since the two transpositions ($\|_1 \leftrightsquigarrow \|_2$) and ($\|_2 \leftrightsquigarrow \|_3$) generate $S_3$, it is sufficient to consider them. In what follows we focus on the transposition ($\|_1 \leftrightsquigarrow \|_2$). For the transposition ($\|_2 \leftrightsquigarrow \|_3$) we just state the result, the proof being identical.

Let us denote by $S(\psi)$ the (orientation reversing) symmetric matrix associated to a reflection of a vector $v \in \mathbb{R}^2 \times \{0\}$ with respect to a horizontal line of $\mathbb{R}^2 \times \{0\}$ which passes through the origin and make an angle $\psi$ with the $e_1$ axis. Let us also denote by $S_\star(\psi)$ the orientation preserving counterpart. In coordinates

$$S(\psi) := \begin{pmatrix} \cos(2\psi) & \sin(2\psi) & 0 \\ \sin(2\psi) & -\cos(2\psi) & 0 \\ 0 & 0 & 1 \end{pmatrix}, \quad S_\star(\psi) := \begin{pmatrix} \cos(2\psi) & \sin(2\psi) & 0 \\ \sin(2\psi) & -\cos(2\psi) & 0 \\ 0 & 0 & -1 \end{pmatrix}. \tag{15}$$

Notice that for every $(\theta, \psi) \in \mathbb{R}^2$ one has $S_\star(\psi) R(\theta) = S_\star(\psi - \theta/2)$, $R(\theta) S_\star(\psi) = S_\star(\psi + \theta/2)$ and

$$R(\theta) S_\star(\psi) R^\mathsf{T}(\theta) = R(\theta) S_\star(\psi + \theta/2) = S_\star(\psi + \theta). \tag{16}$$

Next, let us set $L := [e_2 | e_1 | e_3]$. Note that when the reflection $L$ is applied to the geometric domain described by a swimmer in the initial angular state $\theta_0 = -\pi/6$, the result (up to a translation) is just a swap of the arms $\|_1$ and $\|_2$ of sPr₃, the arm $\|_3$ remaining (still up to a translation) the same (cf. **Figure** 2). Stokes equations then justify the following

CONDITION 4. (SWAP ($\|_1 \leftrightsquigarrow \|_2$)) *Let the initial position be $p_0 := (x_0, y_0, -\pi/6)$. If $\gamma(p_0, L\xi)$ is a solution of the control system (5), then so is $\gamma(S(\pi/2) p_0, \xi)$ and the following relation holds*

$$\gamma(p_0, L\xi) = S_\star\left(\tfrac{\pi}{2}\right) \gamma\left(S\left(\tfrac{\pi}{2}\right) p_0, \xi\right) - \tfrac{\pi}{3} e_3. \tag{17}$$

**Remark 5.** Physically speaking, the previous hypothesis is a consequence of the invariance of Stokes equations with respect to the observation point (cf. **Figure** 2). Indeed, an observer watching the dynamics $\gamma(p_0, \xi)$ of sPr$_3$ projected on a glass, imposes to another observer, lying on the other side of the glass, to watch the dynamics $\gamma(p_0, L\xi)$ of a micro-swimmer obtained from sPr$_3$ by inverting arms $\|_1$ and $\|_2$.

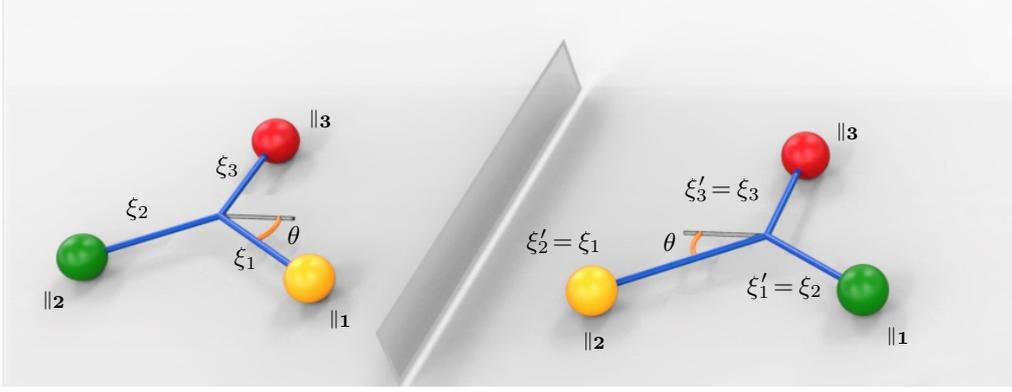

**Figure 2.** An observer watching the dynamics $\gamma(p_0, \xi)$ of sPr$_3$ projected on a glass, imposes to another observer, lying on the other side of the glass, to watch the dynamics $\gamma(p_0, L\xi)$ of a micro-swimmer obtained from sPr$_3$ by inverting arms $\|_1$ and $\|_2$.

**Proposition 6.** *If the control system* (5) *is invariant under the swap* ($\|_1 \rightsquigarrow \|_2$) *and* $T_\xi \mathcal{M} \simeq \mathbb{R}^3$ *for all* $\xi \in \mathcal{M}$, *then*

$$F(L\xi) = S_L \, F(\xi) \, L \quad \forall \xi \in \mathcal{M}, \tag{18}$$

*with* $S_L := S_\star(2\pi/3)$.

**Proof.** Let $\gamma(p_0, L\xi)$ be a solution of the control system (14). From the hypothesis of rotational invariance we get (cf. (10)) that for any $\theta_1 \in \mathbb{R}$

$$\gamma(p_0, L\xi) = R^\mathsf{T}(\theta_1)\gamma(p_0 + \theta_1 e_3, L\xi) - \theta_1 e_3 + A(-\theta_1)p_0. \tag{19}$$

Let us set $\theta_0 := \gamma(p_0, L\xi)(0) \cdot e_3$. By choosing $\theta_1 := -(\theta_0 + \pi/6)$ and setting $p_1 := p_0 + \theta_1 e_3$, we get from Condition 4

$$\gamma(p_0, L\xi) = R^\mathsf{T}(\theta_1) S_\star(\pi/2) \gamma(S(\pi/2)p_1, \xi) - \left(\frac{\pi}{3} + \theta_1\right)e_3 + A(-\theta_1)p_0. \tag{20}$$

Since both $\gamma(p_0, L\xi)$ and $\gamma(S(\pi/2)p_1, \xi)$ are solutions of the control system (14), we get on the one hand

$$\dot\gamma(p_0, L\xi) \stackrel{(14)}{=} R(\gamma(p_0, L\xi) \cdot e_3) F(L\xi) L \dot\xi, \tag{21}$$

and on the other hand, by first using (20) and then (14)

$$\dot\gamma(p_0, L\xi) \stackrel{(20)}{=} R^\mathsf{T}(\theta_1) S_\star(\pi/2) \dot\gamma(S(\pi/2)p_1, \xi)$$
$$\stackrel{(14)}{=} R^\mathsf{T}(\theta_1) S_\star(\pi/2) R(\gamma(S(\pi/2)p_1, \xi) \cdot e_3) F(\xi) \dot\xi. \tag{22}$$

By imposing the equality of (21) and (22) at time $t = 0$, and recalling that $T_{\xi_0}\mathcal{M} \simeq \mathbb{R}^3$, we obtain that $F(L\xi_0) = R^\mathsf{T}(\theta_0 + \theta_1) S_\star(\pi/2) R(\theta_0 + \theta_1) F(\xi_0) L$ for every $(\theta, \xi) \in \mathbb{R} \times \mathcal{M}$. Eventually, by evaluating the previous expression for $\theta_0 + \theta_1 = -\pi/6$, taking into account (16) and the arbitrariness of $\xi_0 \in \mathcal{M}$ we conclude. $\square$

For what concerns the transposition ($\|_2 \rightsquigarrow \|_3$), let us set $M := [e_1 | e_3 | e_2]$. For symmetry reasons we can state the following

CONDITION 7. (SWAP ($\|_\mathbf{2} \leftrightsquigarrow \|_\mathbf{3}$)) *Let the initial position be $p_0 := (x_0, y_0, \pi/2)$. If $\gamma(p_0, M\xi)$ is a solution of the control system* (5)*, then so is $\gamma(S(\pi/2)p_0, \xi)$ and the following relation holds:*

$$\gamma(p_0, M\xi) = S_\star(\pi/2)\gamma(S(\pi/2)p_0, \xi) + \pi e_3. \tag{23}$$

PROPOSITION 8. *If the control system* (5) *is invariant under the swap* ($\|_\mathbf{2} \leftrightsquigarrow \|_\mathbf{3}$) *and $T_\xi \mathcal{M} \simeq \mathbb{R}^3$ for all $\xi \in \mathcal{M}$, then*

$$F(M\xi) = S_M F(\xi) M \quad \forall \xi \in \mathcal{M}, \tag{24}$$

*with $S_M := S_\star(0)$.*

## 4. THE CONTROL SYSTEM IN THE RANGE OF SMALL STROKES

Let us restart from (14). The response of the control system is governed by the matrix valued function $F: \mathbb{R} \times \mathcal{M} \to \mathbb{R}^{3 \times 3}$ which, due to Proposition 3, can be factorised as:

$$F(\theta, \zeta) = R(\theta) F(\zeta) \quad \text{with} \quad F(\zeta) := F(0, \zeta). \tag{25}$$

In what follows we suppose that $\zeta := \xi_0 + \xi$ with $\xi_0 \in \mathcal{M}$ having all its components equal, and we set $F_{\xi_0}(\xi) := F(\xi_0 + \xi)$. Since $F$ is an analytic function (cf. [ADH+13]) we can write the first order expansion[2]

$$F_{\xi_0}(\xi)\eta = F_0 \eta + \mathcal{H}_0(\xi \otimes \eta) + o(|\xi|)\eta, \tag{26}$$

with $F_0 := F(\xi_0) \in \mathbb{R}^{3 \times 3}$ and $\mathcal{H}_0 \in \mathcal{L}(\mathbb{R}^3 \otimes \mathbb{R}^3, \mathbb{R}^3)$ representing the first order derivative of $F_{\xi_0}$ at $\xi = 0$.

The aim of this section is to reveal the structure of the zeroth and first order terms of the expansion (26), in view of the symmetry conditions that $F_{\xi_0}$ must satisfy due to Propositions 6 and 8, i.e:

$$F_{\xi_0}(L\xi) = S_L F(\xi) L \quad, \quad F_{\xi_0}(M\xi) = S_M F_{\xi_0}(\xi) M \quad \forall \xi \in \mathcal{M}. \tag{27}$$

Let us prove the following

LEMMA 9. *Let us suppose that for some matrices $A$ and $S_A$ we have $F_{\xi_0}(A\xi) = S_A F_{\xi_0}(\xi) A$ for every $\xi \in \mathcal{M}$. Then necessarily*

$$F_0 = S_A F_0 A, \tag{28}$$

*and*

$$\mathcal{H}_0((A\xi) \otimes \eta) = S_A \mathcal{H}_0(\xi \otimes (A\eta)) \quad \forall \xi, \eta \in \mathbb{R}^3. \tag{29}$$

**Proof.** Evaluating the condition $F_{\xi_0}(A\xi) = S_A F_{\xi_0}(\xi) A$ at $\xi = 0$, we get (28). Next, By setting $\eta := A\eta$ into the expansion (26) we get

$$F_{\xi_0}(\xi) A\eta = F_0 A\eta + \mathcal{H}_0(\xi \otimes (A\eta)) + o(|\xi|)\eta \quad \forall \xi, \eta \in \mathbb{R}^3. \tag{30}$$

Therefore from (26) and (30) we get

$$\begin{aligned}
\mathcal{H}_0((A\xi) \otimes \eta) &\stackrel{(26)}{=} F_{\xi_0}(A\xi)\eta - F_0 \eta + o(|\xi|)\eta \\
&= S_A F_{\xi_0}(\xi) A\eta - F_0 \eta + o(|\xi|)\eta \\
&\stackrel{(30)}{=} S_A[F_0 A\eta + \mathcal{H}_0(\xi \otimes (A\eta))] - F_0 \eta + o(|\xi|)\eta \\
&\stackrel{(28)}{=} S_A \mathcal{H}_0(\xi \otimes (A\eta)) + o(|\xi|)\eta,
\end{aligned} \tag{31}$$

---

2. Here, for notational convenience, we make use of the universal factorization property of tensor spaces in terms of multilinear maps (cf. [Sch75]).

and hence (29). □

### 4.1. The zeroth order term $F_0$

Applying Lemma 9 to the two matrix actions $A := L$ and $A := M$ (cf. (27)) we respectively get the two conditions $F_0 = S_L F_0 L$ and $F_0 = S_M F_0 M$, which constitute a linear system of matrix equations in the unknown $F_0$. A direct computation shows that the space of solutions is one-dimensional, and has the following structure:

$$F_0 := \begin{pmatrix} -2\mathfrak{a} & \mathfrak{a} & \mathfrak{a} \\ 0 & \sqrt{3}\mathfrak{a} & -\sqrt{3}\mathfrak{a} \\ 0 & 0 & 0 \end{pmatrix} \quad \text{with} \quad \mathfrak{a} \in \mathbb{R}. \tag{32}$$

In connection with optimality questions, it is convenient to introduce the following orthogonal basis of $\mathbb{R}^3$:

$$\tau_1 := (0, -1, 1) \quad , \quad \tau_2 := \frac{1}{\sqrt{3}}(-2, 1, 1) \quad , \quad \tau_3 := (1, 1, 1). \tag{33}$$

The matrix $F_0$ then reads as $F_0 = \mathfrak{a}\sqrt{3}\,[\tau_2|-\tau_1|0]^\mathsf{T}$.

### 4.2. The first order term $\mathcal{H}_0$

Evaluating (29) on the basis $(e_i \otimes e_j)_{i,j \in \mathbb{N}_3}$ we get $\mathcal{H}_0((Ae_i) \otimes e_j) = S_A \mathcal{H}_0(e_i \otimes (Ae_j))$ for every $i$, $j \in \mathbb{N}_3$. In particular, if the action of $A$ on the ordered basis $(e_1, e_2, e_3)$ consists in a permutation of the ordered basis, i.e. if $Ae_i = e_{\sigma(i)}$ for every $i \in \mathbb{N}_3$ and a suitable permutations $\sigma \in \mathbb{N}_3 \to \mathbb{N}_3$, then the structural condition (29) in Lemma 9 reads as

$$\mathcal{H}_0(e_{\sigma(i)} \otimes e_j) = S_A \mathcal{H}_0(e_i \otimes e_{\sigma(j)}) \quad \forall i, j \in \mathbb{N}_3. \tag{34}$$

In particular we have:

COROLLARY 10. *Let us denote by $\sigma_L, \sigma_M \in S_3$ the permutations $\sigma_L = (2, 1, 3)$ and $\sigma_M = (1, 3, 2)$. The following relations hold:*

$$\begin{align}
\mathcal{H}_0(e_{\sigma_L(i)} \otimes e_j) &= S_L \mathcal{H}_0(e_i \otimes e_{\sigma_L(j)}) \tag{35} \\
\mathcal{H}_0(e_{\sigma_M(i)} \otimes e_j) &= S_M \mathcal{H}_0(e_i \otimes e_{\sigma_M(j)}) \tag{36}
\end{align}$$

*for every $i, j \in \mathbb{N}_3$.*

**Proof.** It is sufficient to recall that $F_{\xi_0}$ must satisfy relations (27) and that $Le_i = e_{\sigma_L(i)}$, $Me_i = e_{\sigma_M(i)}$ for every $i \in \mathbb{N}_3$. □

To reveal the structure of $\mathcal{H}_0$ it is necessary to compute the solution space of the system of 18 vector equations (i.e. of 54 scalar equations) given by (35)-(36). This laborious task can be addressed via a symbolic mathematical computation program or, for the brave, with the aid of some clever observation. The first thing to observe is that what we are really interested in, is the structure of the matrices $A_k$ defined for any $k \in \mathbb{N}_3$ by the position $A_k := (\mathcal{H}_0(e_i \otimes e_j) \cdot e_k)_{i,j \in \mathbb{N}_3}$. Indeed, for any $\xi, \eta \in \mathbb{R}^3$ the vector $\mathcal{H}_0(\xi \otimes \eta) \in \mathbb{R}^3$ is given by

$$\mathcal{H}_0(\xi \otimes \eta) = \sum_{k \in \mathbb{N}_3} (A_k \eta \cdot \xi)\, e_k. \tag{37}$$

Next we observe that since $S_L$ and $S_M$ are idempotent matrices multiplying both members of (35) and (36) by $S_L$ and $S_M$ we respectively get the relations $\mathcal{H}_0(e_i \otimes e_{\sigma_L(j)}) = S_L \mathcal{H}_0(e_{\sigma_L(i)} \otimes e_j)$ and $\mathcal{H}_0(e_i \otimes e_{\sigma_M(j)}) = S_M \mathcal{H}_0(e_{\sigma_M(i)} \otimes e_j)$ for any $i, j \in \mathbb{N}_3$. This reduces the set of vector equations from 18 to 10. Eventually, with the help of relations (16), one discovers that by setting

$$\alpha := \mathcal{H}_0(e_1 \otimes e_2) \cdot e_1 + \frac{1}{2}\mathcal{H}_0(e_3 \otimes e_2) \cdot e_1 \quad \text{and} \quad \gamma := \mathcal{H}_0(e_1 \otimes e_2) \cdot e_3, \tag{38}$$

as well as $\beta := \frac{3}{2}\mathcal{H}_0(e_3 \otimes e_2) \cdot e_1$ and $\lambda := \mathcal{H}_0(e_1 \otimes e_1) \cdot e_1$, one gets

$$A_1 = \begin{pmatrix} \lambda & \alpha - \frac{1}{3}\beta & \alpha - \frac{1}{3}\beta \\ -\alpha - \frac{1}{3}\beta & -\frac{\lambda}{2} & \frac{2}{3}\beta \\ -\alpha - \frac{1}{3}\beta & \frac{2}{3}\beta & -\frac{\lambda}{2} \end{pmatrix}, \quad A_2 = \sqrt{3}\begin{pmatrix} 0 & \frac{\alpha+\beta}{3} & -\frac{\alpha+\beta}{3} \\ \frac{\beta-\alpha}{3} & -\frac{\lambda}{2} & -\frac{2\alpha}{3} \\ \frac{\alpha-\beta}{3} & \frac{2\alpha}{3} & \frac{\lambda}{2} \end{pmatrix}, \tag{39}$$

and

$$A_3 = \begin{pmatrix} 0 & \gamma & -\gamma \\ -\gamma & 0 & \gamma \\ \gamma & -\gamma & 0 \end{pmatrix}. \tag{40}$$

In particular, the skew-symmetric parts (denoted by $M_1, M_2, M_3$) of the matrices $A_1, A_2, A_3$ which, as we shall see, are the only ones to contribute to the net displacement of the micro-swimmer, are given by

$$M_1 := \begin{pmatrix} 0 & \alpha & \alpha \\ -\alpha & 0 & 0 \\ -\alpha & 0 & 0 \end{pmatrix}, \quad M_2 := \frac{1}{\sqrt{3}}\begin{pmatrix} 0 & \alpha & -\alpha \\ -\alpha & 0 & -2\alpha \\ \alpha & 2\alpha & 0 \end{pmatrix}, \quad M_3 := A_3. \tag{41}$$

Let us note that the orthogonal basis (33) is orientation-reversing and that the matrices in (41) can be characterized by the actions

$$M_1 \xi = \alpha \xi \times \tau_1, \quad M_2 \xi = \alpha \xi \times \tau_2, \quad M_3 \xi = \gamma \xi \times \tau_3 \quad \forall \xi \in \mathbb{R}^3. \tag{42}$$

### 4.3. The linearized control equations

In what follows we denote by $I$ the closed interval $[0, 2\pi] \subset \mathbb{R}$ and by $H^1_\sharp(I, \mathbb{R}^3)$ the so-called **strokes space**, mathematically identified with the Sobolev space of $2\pi$-periodic vector valued functions of $L^2_\sharp(I, \mathbb{R}^3)$ having first order weak derivative in $L^2_\sharp(I, \mathbb{R}^3)$. For every $p \in L^1_\sharp(I, \mathbb{R}^3)$ we denote by $\langle p \rangle := (2\pi)^{-1} \int_I p(s) \, \mathrm{d}s$ the average of $p$ on $I$.

In section 3 we have shown that the control system governing the evolution of the swimmer sPr₃ under the action of the control parameters $\zeta \in \mathcal{M}$ can be written as $\dot{p} = R(\theta)F(\zeta)\dot{\zeta}$, where the system $F: \mathcal{M} \to \mathbb{R}^{3\times 3}$ is given by (25) and $\dot{\zeta} \in T_\zeta \mathcal{M}$. Moreover, cf. (26) and (37), if we set $\zeta = \xi_0 + \xi$, the behaviour of the system around $\xi = 0$, up to higher order terms, is given by

$$\dot{p} = R(\theta) F_0 \dot{\xi} + R(\theta) \sum_{k \in \mathbb{N}_3} (A_k \dot{\xi} \cdot \xi) e_k. \tag{43}$$

In particular, denoting by $(c, \theta) \in \mathbb{R}^2 \times \mathbb{R}$ the components of $p$ and taking into account that $F_0^\mathsf{T} e_3 = 0$ because of (32), we get

$$\dot{\theta} = A_3 \dot{\xi} \cdot \xi = M_3 \dot{\xi} \cdot \xi = \gamma \tau_3 \cdot (\xi \times \dot{\xi}), \tag{44}$$

which can be easily integrated:

$$\theta(t) - \theta_0 = \tau_3 \cdot \gamma \int_0^t \xi(s) \times \dot{\xi}(s) \, \mathrm{d}s \quad \text{with} \quad \theta_0 := \theta(0).$$

Hence the net angular displacement $\delta\theta$ corresponding to a unit stroke $\xi$ is given by

$$\delta\theta := 2\pi \langle \dot{\theta} \rangle = 2\gamma\pi \langle \xi \times \dot{\xi} \rangle \cdot \tau_3. \tag{45}$$

Since the dynamics of $\theta$ does not depend on the one of $c$, it is convenient to consider the $2d$ projections of the matrices $R(\theta)$ and $F_0^\mathsf{T}$. In what follows we denote by $\hat{e}_1, \hat{e}_2$ the standard basis of $\mathbb{R}^2$ and to lighten notation, since no confusion may arise, we still denote by $R(\theta)$ and $F_0$ the $2d$ projections

$$F_0 := \mathfrak{a} \begin{pmatrix} -2 & 1 & 1 \\ 0 & \sqrt{3} & -\sqrt{3} \end{pmatrix}, \quad R(\theta) := \begin{pmatrix} \cos\theta & -\sin\theta \\ \sin\theta & \cos\theta \end{pmatrix}, \tag{46}$$

so that the dynamics of $c$ is described by the system

$$\dot{c} = R(\theta) F_0 \dot{\xi} + R(\theta) \sum_{j \in \mathbb{N}_2} (A_j \dot{\xi} \cdot \xi) \hat{e}_j. \tag{47}$$

Next, we observe that the first order expansion of the $2d$ rotation matrix $R(\theta)$ around $\theta = 0$ gives

$$R(\theta) = \text{Id} + K(\theta) + \mathcal{O}(\theta^2) \quad \text{with} \quad K(\theta) := \begin{pmatrix} 0 & -\theta \\ \theta & 0 \end{pmatrix}, \tag{48}$$

and therefore, up to higher order terms in $\theta$, we get

$$\dot{c} = (\text{Id} + K(\theta))(F_0 \dot{\xi} + \sum_{j \in \mathbb{N}_2} (A_j \dot{\xi} \cdot \xi) \hat{e}_j) + \mathcal{O}(\theta^2). \tag{49}$$

An integration on $I = [0, 2\pi]$ of the previous relation gives an estimate of the net displacement $\delta c$ undergone by the center $c$ of sPr$_3$ in correspondence to a small stroke. Moreover, relation (49) allows to interpret $\delta c$ as a map $\delta c : \xi \in H^1_\sharp(I, \mathbb{R}^3) \mapsto 2\pi \langle \dot{c}(\xi) \rangle \in \mathbb{R}^3$.

PROPOSITION 11. *For any $\xi \in H^1_\sharp(I, \mathbb{R}^3)$, $\xi : I \to \mathcal{M}_0$, in a neighbourhood of $0 \in H^1_\sharp(I, \mathbb{R}^3)$ the following estimate holds*

$$\delta c(\xi) = 2\pi \langle A_1 \dot{\xi} \cdot \xi \rangle \hat{e}_1 + 2\pi \langle A_2 \dot{\xi} \cdot \xi \rangle \hat{e}_2 + \mathcal{O}(\|\xi\|_{H^1_\sharp}^3). \tag{50}$$

*In physical terms: in the limit of small strokes around a constant reference shape $\xi_0$ the net displacement is given by* (50).

**Proof.** We first note that the term $\langle F_0 \dot{\xi} \rangle$ is zero by periodicity. Next we observe that is sufficient to prove the estimate for the scalar terms of the form $\langle \theta(\xi, \cdot) \dot{\xi}_i \rangle$ and $\langle \theta(\xi, \cdot) A_j \dot{\xi} \cdot \xi \rangle$ with $i \in \mathbb{N}_3$, $j \in \mathbb{N}_2$ and

$$\theta(\xi, t) := \int_0^t M_3 \dot{\xi}(s) \cdot \xi(s) \, ds \quad \forall t \in I. \tag{51}$$

Let us focus on the terms of the form $\langle \theta(\xi, \cdot) \dot{\xi}_i \rangle$, the estimate concerning the other ones can be treated in exactly the same way. We have

$$\begin{aligned}
\int_I \theta(\xi, t) \dot{\xi}_i \, dt &= \int_I \left( \int_0^t M_3 \dot{\xi}(s) \cdot \xi(s) \, ds \right) \dot{\xi}_i(t) \, dt \\
&= \int_0^{2\pi} M_3 \dot{\xi}(t) \cdot \xi(t) \int_t^{2\pi} \dot{\xi}_i(s) \, ds \, dt \\
&\leqslant \gamma |\tau_3| \int_I |\dot{\xi}(t)| \cdot |\xi(t) - \xi(0)|^2 \, dt.
\end{aligned} \tag{52}$$

The Sobolev-Morrey embedding $H^1_\sharp(I, \mathbb{R}^3) \subseteq L^\infty_\sharp(I, \mathbb{R}^3)$ gives the existence of a $c_S > 0$ such that $\|\xi\|_\infty \leqslant c_S \|\xi\|_{H^1_\sharp}$ for every $\xi \in H^1_\sharp(I, \mathbb{R}^3)$. Therefore for some $c_\gamma > 0$ depending on $M_3$ only, we get

$$\begin{aligned}
\int_I \theta(\xi, t) \dot{\xi}_i \, dt &\leqslant c_\gamma \|\xi\|_\infty^2 \|\xi\|_{H^1_\sharp(I, \mathbb{R}^3)} \\
&\leqslant c_\gamma c_S \|\xi\|_{H^1_\sharp(I, \mathbb{R}^3)}^3.
\end{aligned} \tag{53}$$

The proof is complete. $\square$

Collecting (44) and (50) we may therefore assume the net displacement $\delta p := p(2\pi) - p(0)$, undergone by the position $p$ of sPr$_3$ in correspondence to a small stroke $\xi$, to be given by (cf. (42))

$$\begin{aligned}
\frac{1}{2\pi} \delta p &= \langle M_1 \dot{\xi} \cdot \xi \rangle e_1 + \langle M_2 \dot{\xi} \cdot \xi \rangle e_2 + \langle M_3 \dot{\xi} \cdot \xi \rangle e_3 \\
&= \alpha [\langle \xi \times \dot{\xi} \rangle \cdot \tau_1] e_1 + \alpha [\langle \xi \times \dot{\xi} \rangle \cdot \tau_2] e_2 + \gamma [\langle \xi \times \dot{\xi} \rangle \cdot \tau_3] e_3,
\end{aligned} \tag{54}$$

where in writing the previous expression we have taken into account that only the skew-symmetric part of the matrices $A_k$ contribute to the displacement. Indeed, if we denote by $A_k^{\text{sym}}$ the symmetric part of the matrix $A_k$, integrating by parts we get $\langle A_k^{\text{sym}} \xi \cdot \dot{\xi} \rangle = \langle A_k^{\text{sym}} \dot{\xi} \cdot \xi \rangle = -\langle A_k^{\text{sym}} \xi \cdot \dot{\xi} \rangle$ and therefore $\langle A_k^{\text{sym}} \dot{\xi} \cdot \xi \rangle = 0$ for any $\xi \in H^1_\sharp(I, \mathbb{R}^3)$.

## 5. ENERGY MINIMIZING STROKES

Following the swimming efficiency suggested by Lighthill [Lig52], we adopt the following notion of optimality: energy minimizing strokes are the ones that minimize the kinematic energy dissipated while trying to reach a given net displacement $\delta p \in \mathbb{R}^3$ in one stroke. In mathematical terms, the total energy dissipation due to a stroke $\xi \in H^1_\sharp(I, \mathbb{R}^3)$, $\xi\colon I \to \mathcal{M}$, can be evaluated by the means of a suitable quadratic energy functional

$$\mathcal{G}(\xi) := \int_I \mathfrak{g}(\xi(t)) \dot\xi(t) \cdot \dot\xi(t) \,\mathrm{d}t, \tag{55}$$

in which the energy density $\mathfrak{g} \in C^1(\mathbb{R}^3)$ is a function taking values in the space of symmetric and definite positive matrices of $\mathbb{R}^{3\times 3}$ (cfr. [DAL12]). In the limit of small strokes one considers the energy $\mathcal{G}$ as arising from the approximation $\mathfrak{g}(\xi) = \mathfrak{g}(0) + o(1)$, in which $\mathfrak{g}(0) \in \mathbb{R}^{3\times 3}$ is a symmetric and definite positive matrix. Namely

$$\mathcal{G}(\xi) := \int_I Q_\mathfrak{g}(\dot\xi(t)) \,\mathrm{d}t \tag{56}$$

where we have denoted by $Q_\mathfrak{g}(\eta)$ the definite positive quadratic form given by $Q_\mathfrak{g}(\eta) := \mathfrak{g}(0)\eta \cdot \eta$.

For the same symmetry reasons discussed in section 3, for every $\eta \in \mathbb{R}^3$ the function $Q_\mathfrak{g}$ must satisfy the relations

$$Q_\mathfrak{g}(L\eta) = Q_\mathfrak{g}(M\eta) = Q_\mathfrak{g}(\eta), \tag{57}$$

in which $L = [e_2|e_1|e_3]$ and $M = [e_1|e_3|e_2]$. A straightforward computation shows that the previous conditions imply the existence of two parameters $h$ and $\kappa > \max(h, -2h)$, such that the symmetric and definite positive matrix $G$ which represents $Q_\mathfrak{g}$ is given by

$$G = \begin{pmatrix} \kappa & h & h \\ h & \kappa & h \\ h & h & \kappa \end{pmatrix}. \tag{58}$$

Let us note that one has $G\tau_1 = (\kappa - h)\tau_1$, $G\tau_2 = (\kappa - h)\tau_2$ and $G\tau_3 = (\kappa + 2h)\tau_3$, so that, up to a rescaling, the orthogonal basis $(\tau_1, \tau_2, \tau_3)$ is invariant and orthogonal with respect to $G$. It is convenient to denote by $\mathfrak{g}_1 := \mathfrak{g}_2 := (\kappa - h)$ and $\mathfrak{g}_3 := (\kappa + 2h)$ the eigenvalues of $G$. Since $G$ is symmetric and definite positive, the coefficients $(\mathfrak{g}_1, \mathfrak{g}_2, \mathfrak{g}_3)$ can be interpreted as the metric coefficients of the inner product defined for any $a, b \in \mathbb{R}^3$ by $(a,b)_\mathfrak{g} := 2\pi \Lambda_\mathfrak{g} a \cdot b$, with $\Lambda_\mathfrak{g} := \mathrm{diag}(\mathfrak{g}_i)$. Let us observe that the metric coefficients $(\mathfrak{g}_i)_{i \in \mathbb{N}_3}$ permit to diagonalize $G$ in the form

$$G = U\Lambda_\mathfrak{g} U^\mathsf{T}, \quad U := [\hat\tau_1|\hat\tau_2|\hat\tau_3], \quad \Lambda_\mathfrak{g} := \mathrm{diag}(\mathfrak{g}_i). \tag{59}$$

with $(\hat\tau_1, \hat\tau_2, \hat\tau_3)$ being the normalization of $(\tau_1, \tau_2, \tau_3)$ with respect to the usual Euclidean metric. It is worth noting that the basis $(\hat\tau_1, \hat\tau_2, \hat\tau_3)$ is both orthogonal and $\mathfrak{g}$-orthogonal.

After that we can write the following equivalent expression of the functional $\mathcal{G}$:

$$\mathcal{G}\colon \xi \in H^1_\sharp(I, \mathbb{R}^3) \mapsto \sum_{i \in \mathbb{N}_3} \mathfrak{g}_i \int_I (\dot\xi(t) \cdot \hat\tau_i)^2 \,\mathrm{d}t. \tag{60}$$

We aim at minimizing $\mathcal{G}$ in $H^1_\sharp(I, \mathbb{R}^3)$ subject to a prescribed net displacement $\delta p \in \mathbb{R}^3$, i.e. (cf. (53)) subject to the constraint

$$\sum_{i \in \mathbb{N}_3} \left( \mathfrak{h}_i \hat\tau_i \cdot \int_I \xi(t) \times \dot\xi(t) \,\mathrm{d}t \right) e_i = \delta p, \tag{61}$$

with $\mathfrak{h}_1/|\tau_1| = \mathfrak{h}_2/|\tau_2| = \alpha$ and $\mathfrak{h}_3/|\tau_3| = \gamma$. More precisely, we are going to prove the following

THEOREM 12. *Let $\delta p \in \mathbb{R}^3$ be a prescribed displacement. Any minimizer $\xi \in H^1_\sharp(I, \mathbb{R}^3)$ of the energy functional* (60) *subject to the constraint* (61) *is of the type*

$$\xi(t) := (\cos t) u + (\sin t) v, \tag{62}$$

*i.e. an ellipse of $\mathbb{R}^3$ centered at the origin and contained in the plane spanned by the vectors $u$ and $v$. The vectors $u, v \in \mathbb{R}^3$ are determined as follows:*

- *We compute the vector $\omega$ via the relation*

$$\omega := \mathrm{diag}\left(\frac{\sqrt{\mathfrak{g}\mathfrak{g}_\theta}}{\mathfrak{h}}, \frac{\sqrt{\mathfrak{g}\mathfrak{g}_\theta}}{\mathfrak{h}}, \frac{\mathfrak{g}}{\mathfrak{h}_\theta}\right)\delta p, \qquad (63)$$

*in which we have made use of the notation $\mathfrak{g}_1 = \mathfrak{g}_2 = \mathfrak{g}$, $\mathfrak{g}_3 = \mathfrak{g}_\theta$, and similarly $\mathfrak{h}_1 = \mathfrak{h}_2 = \mathfrak{h}$, $\mathfrak{h}_3 = \mathfrak{h}_\theta$. We then consider a vector $\sigma \in \mathbb{R}^3$ in the plane orthogonal to $\omega$ and such that*

$$|\sigma|^2 = |\omega|, \qquad (64)$$

*e.g. $\sigma := \sqrt{|\omega|}\frac{\mu \times \omega}{|\mu \times \omega|}$, with $\mu$ linearly independent from $\omega$.*

- *We set $\hat{\omega} := \omega/|\omega|$ and compute the vectors $u$ and $v$ via the relations*

$$u := \frac{U\Lambda_\mathfrak{g}^{-1/2}}{\sqrt{2\pi}}\sigma \quad , \quad v := \frac{U\Lambda_\mathfrak{g}^{-1/2}}{\sqrt{2\pi}}(\sigma \times \hat{\omega}). \qquad (65)$$

*We then have $v \times u = \omega$ and the minimum value of $\mathcal{G}$ is equal to $|\omega|$.*

*Moreover, the vectors $u$ and $v$ are $\mathfrak{g}$-orthogonal, i.e. orthogonal with respect to the inner product defined for every $a, b \in \mathbb{R}^3$ by $(a, b)_\mathfrak{g} := 2\pi\Lambda_\mathfrak{g} a \cdot b$, and have the same $\mathfrak{g}$-norm: $|u|_\mathfrak{g}^2 = |v|_\mathfrak{g}^2 = |\omega|$.*

**Remark 13.** If $\delta p \parallel e_i$, with $i \in \mathbb{N}_3$, then $\hat{\omega} = e_i$. Hence, setting $\sigma = \sqrt{|\omega|}e_{i-1}$ relation (64) is obviously satisfied and from (65) we get (the notation make use of cyclic permutations of the indices):

$$u := \sqrt{\frac{|\omega|}{2\pi\mathfrak{g}_{i+1}}}\hat{\tau}_{i-1} \quad , \quad v := \sqrt{\frac{|\omega|}{2\pi\mathfrak{g}_{i-1}}}\hat{\tau}_{i+1}. \qquad (66)$$

Thus, an energy minimizing net displacement along the $x$-asxis direction is achieved, with respect to the standard euclidean space $(\mathbb{R}^3, (\cdot, \cdot)_2)$, via an elliptic stroke contained in the plane orthogonal to the vector $\hat{\tau}_1$. Similarly a pure along $y$ (resp. along $\theta$) net displacement is achieved via an elliptic stroke contained in the plane orthogonal to $\hat{\tau}_2$ (resp. to $\hat{\tau}_3$). On the other hand, with respect to the inner-product space $(\mathbb{R}^3, (\cdot, \cdot)_\mathfrak{g})$, the energy minimizing strokes describe circles of radius $\sqrt{|\omega|}$ (cf. **Figure** 3).

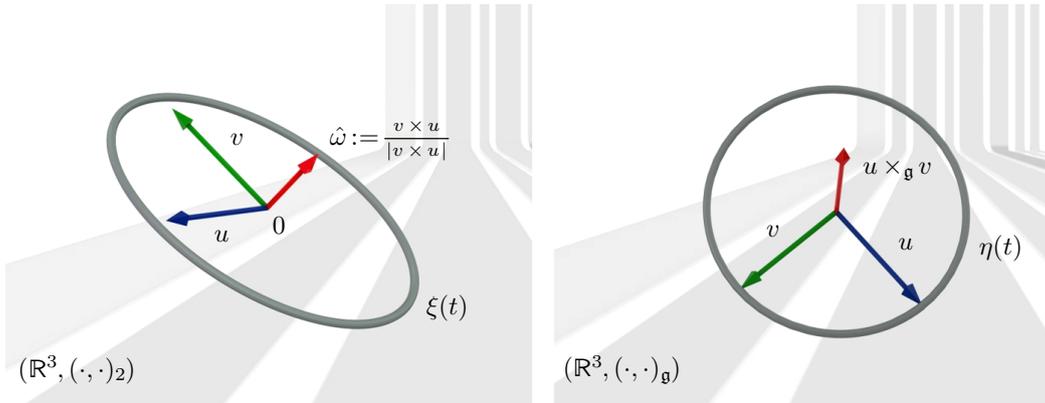

**Figure 3.** (left) With respect to the standard euclidean space $(\mathbb{R}^3, (\cdot, \cdot)_2)$, the energy minimizing strokes able to reach a prescribed displacement $\delta p$ are ellipses of $\mathbb{R}^3$ centered at the origin and contained in the plane spanned by the vectors $u := Ua_1$ and $v := Ub_1$. (right) With respect to the inner-product space $(\mathbb{R}^3, (\cdot, \cdot)_\mathfrak{g})$, the energy minimizing strokes describe circles of radius $\sqrt{|\omega|}$.

The strategy of the proof of Theorem 12 consists in reducing the $H^1_\sharp(I,\mathbb{R}^3)$ constrained minimization problem to a finite dimensional one. To this end we divide the proof in several steps, each step being developed in a different subsection. For convenience of the reader, we briefly summarize the main steps. In STEP$_1$, developed in section 5.1, we rewrite the functional $\mathcal{G}$ and the constraint in terms of the linear coordinate system induced by the metric matrix $G$. We denote this new functional by $\mathcal{G}_U$. In STEP$_2$, developed in section 5.2, we recast the minimization problem for $\mathcal{G}_U$ in the framework of Fourier series. This permits to pass from a minimization problem in $H^1_\sharp(I,\mathbb{R}^3)$ to a minimization problem for a functional $\mathcal{F}$ defined in $\ell^2(\mathbb{R}^3)\times\ell^2(\mathbb{R}^3)$. In STEP$_3$, developed in sections 5.3 and 5.4, we reconduct the minimization problem for $\mathcal{F}$ to a minimization problem for a function $f$ defined in $\mathbb{R}^3\times\mathbb{R}^3$. We then characterize the (global) minimizers of $f$. Eventually, in STEP$_4$, developed in section 5.5, we complete the proof of Theorem 12.

### 5.1. G-orthogonalization

We start by rewriting the functional (60) and the constraint (61) in the linear coordinate system induced by $G$. By the change of variables $\eta(t):=U^\mathsf{T}\xi(t)\in H^1_\sharp(I,\mathbb{R}^3)$, the energy functional $\mathcal{G}$ can be read as

$$\mathcal{G}_U(\eta) = \int_I \Lambda_\mathfrak{g}\dot\eta(t)\cdot\dot\eta(t)\,\mathrm{d}t. \tag{67}$$

with $\mathcal{G}_U(\eta):=\mathcal{G}(\xi)=\mathcal{G}(U\eta)$. For what concerns the constraint (61), it reads as

$$\int_I \dot\eta(t)\times\eta(t)\,\mathrm{d}t = \Lambda_\mathfrak{h}^{-1}\delta p \quad\text{with}\quad \Lambda_\mathfrak{h}:=\mathrm{diag}(\mathfrak{h}_i). \tag{68}$$

Indeed, setting $\eta(t):=U^\mathsf{T}\xi(t)$ into (61) and observing that $\det U = -1$ we get

$$\delta p\cdot e_i = -\mathfrak{h}_i U e_i\cdot U\int_I (U^\mathsf{T}\xi(t))\times(U^\mathsf{T}\dot\xi(t))\,\mathrm{d}t \tag{69}$$

$$= \Lambda_\mathfrak{h}\left(\int_I \dot\eta(t)\times\eta(t)\,\mathrm{d}t\right)\cdot e_i, \tag{70}$$

from which (68) follows.

### 5.2. Fourier approach to the constrained minimization problem: from $H^1_\sharp(I,\mathbb{R}^3)$ to $\ell^2(\mathbb{R}^3)\times\ell^2(\mathbb{R}^3)$

Let us denote by $\ell_2(\mathbb{R}^3)$ the Hilbert space constituted by the sequences $\boldsymbol{u}:=(u_n)_{n\in\mathbb{N}}$ in $\mathbb{R}^3$, such that the norm

$$\|\boldsymbol{u}\|_{\ell_2(\mathbb{R}^3)} := \sum_{n\in\mathbb{N}} |u_n|^2 \tag{71}$$

is finite. We then denote by $\dot\ell_2(\mathbb{R}^3)$ the Hilbert space of sequences $\boldsymbol{u}=(u_n)_{n\in\mathbb{N}}$ of $\ell_2(\mathbb{R}^3)$ such that $(nu_n)_{n\in\mathbb{N}}\in\ell_2(\mathbb{R}^3)$. Every element $\eta\in H^1_\sharp(I,\mathbb{R}^3)$, being $2\pi$-periodic, can be expressed in terms of its Fourier series as

$$\eta(t):=\frac{1}{2}a_0 + \sum_{n\in\mathbb{N}} \cos(nt)a_n + \sin(nt)b_n, \tag{72}$$

with $(a_n,b_n)\in\dot\ell_2(\mathbb{R}^3)\times\dot\ell_2(\mathbb{R}^3)$. Substituting the Fourier series of $\dot\eta$ into the expression (67) of $\mathcal{G}$, we get, thanks to Parseval equality

$$\mathcal{G}_U(\eta) := \int_I \Lambda_\mathfrak{g}\dot\eta(t)\cdot\dot\eta(t)\,\mathrm{d}t = \pi\sum_{n\in\mathbb{N}} n^2(\Lambda_\mathfrak{g}a_n\cdot a_n + \Lambda_\mathfrak{g}b_n\cdot b_n) \tag{73}$$

$$= \tfrac{1}{2}\|\boldsymbol{u}\|^2_{\ell_2(\mathbb{R}^3)} + \tfrac{1}{2}\|\boldsymbol{v}\|^2_{\ell_2(\mathbb{R}^3)}, \tag{74}$$

where we made use of the notation

$$\boldsymbol{u}:=(u_n)_{n\in\mathbb{N}}:=\sqrt{2\pi\Lambda_\mathfrak{g}}(na_n)_{n\in\mathbb{N}} \quad\text{and}\quad \boldsymbol{v}:=(v_n)_{n\in\mathbb{N}}:=\sqrt{2\pi\Lambda_\mathfrak{g}}(nb_n)_{n\in\mathbb{N}}. \tag{75}$$

Let us note that $(\boldsymbol{u}, \boldsymbol{v}) \in \ell_2(\mathbb{R}^3) \times \ell_2(\mathbb{R}^3)$. Next, taking into account the $L^2_\sharp(I)$-orthogonality of the Fourier trigonometric system, it is possible to express the constraint (68) in terms of Fourier coefficients as $2\pi \sum_{n \in \mathbb{N}} \frac{1}{n}(n b_n) \times (n a_n) = \Lambda_\mathfrak{h}^{-1} \delta p$. Moreover, for any $n \in \mathbb{N}$ one has

$$\big(\sqrt{\Lambda_\mathfrak{g}} b_n\big) \times \big(\sqrt{\Lambda_\mathfrak{g}} a_n\big) = \sqrt{(\det \Lambda_\mathfrak{g})\, \Lambda_\mathfrak{g}^{-1}}(b_n \times a_n), \tag{76}$$

and therefore the constraint can ultimately be written as

$$\sum_{n \in \mathbb{N}} \frac{v_n \times u_n}{n} = \sqrt{(\det \Lambda_\mathfrak{g})\, \Lambda_\mathfrak{g}^{-1}} \Lambda_\mathfrak{h}^{-1} \delta p. \tag{77}$$

PROPOSITION 14. *The $H^1_\sharp(I, \mathbb{R}^3)$ minimization of the functional $\mathcal{G}_U$ given by (67) under the constraint (68) is equivalent to the minimization of the functional*

$$\mathcal{F}(\boldsymbol{u}, \boldsymbol{v}) := \frac{1}{2} \|\boldsymbol{u}\|^2_{\ell^2(\mathbb{R}^3)} + \frac{1}{2} \|\boldsymbol{v}\|^2_{\ell^2(\mathbb{R}^3)}, \tag{78}$$

*defined in the product Hilbert space $\ell_2(\mathbb{R}^3) \times \ell_2(\mathbb{R}^3)$ and subject to the constraint*

$$\sum_{n \in \mathbb{N}} \frac{1}{n} v_n \times u_n = \omega \quad \text{with} \quad \omega := \sqrt{(\det \Lambda_\mathfrak{g})\, \Lambda_\mathfrak{g}^{-1}} \Lambda_\mathfrak{h}^{-1} \delta p, \tag{79}$$

*where $\delta p \in \mathbb{R}^3$ is a prescribed net displacement of position.*

## 5.3. From $\ell^2(\mathbb{R}^3) \times \ell^2(\mathbb{R}^3)$ to $\mathbb{R}^3 \times \mathbb{R}^3$

Let us denote by $\boldsymbol{e}_i$ th $i$-th element of the canonical basis of $\ell_2(\mathbb{R})$, i.e. the sequence defined by the position $\boldsymbol{e}_i := (\delta^i_n)_{n \in \mathbb{N}}$, with $\delta^i_n$ Kronecker delta. We have

PROPOSITION 15. *For any $(\boldsymbol{u}, \boldsymbol{v}) \in \ell^2(\mathbb{R}^3) \times \ell^2(\mathbb{R}^3)$ such that the constraint relation (79) holds, there exists two vectors $(u, v) \in \mathbb{R}^3 \times \mathbb{R}^3$ such that once defined the sequences $\boldsymbol{u}_\star := \boldsymbol{e}_1 u$ and $\boldsymbol{v}_\star := \boldsymbol{e}_1 v$ of $\ell^2(\mathbb{R}^3)$, one has*

$$\mathcal{F}(\boldsymbol{u}_\star, \boldsymbol{v}_\star) = \mathcal{F}(\boldsymbol{u}, \boldsymbol{v}) \quad \text{and} \quad v \times u = \omega. \tag{80}$$

**Proof.** If $\omega = 0$ the proof is trivial. Hence, let us denote by $\hat{\omega}$ the unit vector associated to $\omega$. We then consider a couple $(\boldsymbol{u}, \boldsymbol{v}) \in \ell_2(\mathbb{R}^3) \times \ell_2(\mathbb{R}^3)$ and set $\boldsymbol{u}_\star := \boldsymbol{e}_1 u$ and $\boldsymbol{v}_\star := \boldsymbol{e}_1 v$ with $u, v \in \mathbb{R}^3$ chosen so that the following relations hold

$$|u| = \|\boldsymbol{u}\|_{\ell^2(\mathbb{R}^3)} \quad , \quad |v| = \|\boldsymbol{v}\|_{\ell^2(\mathbb{R}^3)} \quad , \quad \frac{u \times v}{|u \times v|} = \hat{\omega}. \tag{81}$$

Since $u \times v = \|\boldsymbol{u}\|_{\ell^2(\mathbb{R}^3)} \|\boldsymbol{v}\|_{\ell^2(\mathbb{R}^3)} (\sin \psi) \hat{\omega}$ for a suitable angle $\psi \in (0, \pi)$, the equality $u \times v = \omega$ can be satisfied by choosing the angle $\psi$ in such a way that

$$\sin \psi = \frac{|\omega|}{\|\boldsymbol{u}\|_{\ell^2(\mathbb{R}^3)} \|\boldsymbol{v}\|_{\ell^2(\mathbb{R}^3)}}, \tag{82}$$

and the previous equation in the variable $\psi$ has a solution as far as the right hand side is less not greater than one. This is indeed the case because due to Cauchy-Schwarz inequality in $\ell^2(\mathbb{R}^3)$ one has

$$|\omega| \leqslant \sum_{n \in \mathbb{N}} \frac{1}{n} |v_n \times u_n| \leqslant \sum_{n \in \mathbb{N}} |v_n \times u_n| \leqslant \|\boldsymbol{u}\|_{\ell^2(\mathbb{R}^3)} \|\boldsymbol{v}\|_{\ell^2(\mathbb{R}^3)}. \tag{83}$$

Eventually, from (81) we have

$$\mathcal{F}(\boldsymbol{u}_\star, \boldsymbol{v}_\star) = \frac{1}{2} |u|^2 + \frac{1}{2} |v|^2 = \frac{1}{2} \|\boldsymbol{u}\|^2_{\ell^2(\mathbb{R}^3)} + \frac{1}{2} \|\boldsymbol{v}\|^2_{\ell^2(\mathbb{R}^3)} = \mathcal{F}(\boldsymbol{u}, \boldsymbol{v}) \tag{84}$$

and this concludes the proof. □

As an immediate consequence of Proposition 15 we get

COROLLARY 16. *The minimization problem for $\mathcal{F}$ in $\ell_2(\mathbb{R}^3) \times \ell_2(\mathbb{R}^3)$, subject to the constraint (79), is equivalent to the minimization in $\mathbb{R}^3 \times \mathbb{R}^3$ of the function*

$$f(u,v) := \frac{1}{2}|u|^2 + \frac{1}{2}|v|^2 \tag{85}$$

*subject to the constraint*

$$v \times u = \omega. \tag{86}$$

**Proof.** It is sufficient to observe that if we denote by $\mathcal{V}_\omega$ the subset of $(\boldsymbol{u}, \boldsymbol{v}) \in \ell_2(\mathbb{R}^3) \times \ell_2(\mathbb{R}^3)$ satisfying the constraint relation (79), and by $V_\omega$ the subset of $(u,v) \in \mathbb{R}^3 \times \mathbb{R}^3$ such that $v \times u = \omega$, then from Proposition 15 we have

$$\min_{(\boldsymbol{u},\boldsymbol{v}) \in \mathcal{V}_\omega} \mathcal{F}(\boldsymbol{u},\boldsymbol{v}) = \min_{(u,v) \in V_\omega} \mathcal{F}(\boldsymbol{e}_1 u, \boldsymbol{e}_1 v) = \min_{(u,v) \in V_\omega} f(u,v) \tag{87}$$

with $\boldsymbol{e}_1 = (\delta_n^1)_{n \in \mathbb{N}}$ and hence $\boldsymbol{e}_1 u = (\delta_n^1 u)_{n \in \mathbb{N}}$. □

### 5.4. The finite dimensional minimization: minimization of $f$

PROPOSITION 17. *Any couple of vectors $(u_\star, v_\star) \in \mathbb{R}^3 \times \mathbb{R}^3$ minimizing the function $f$ given by (85) and subject to the constraint (86), is characterized by the following conditions*

$$|u_\star|^2 = |v_\star|^2 = |\omega| \quad \text{and} \quad u_\star \cdot v_\star = 0. \tag{88}$$

*Therefore, for any $\sigma \in \mathbb{R}^3$ such that $\sigma \cdot \hat{\omega} = 0$ and $|\sigma|^2 = |\omega|$, the couple*

$$(\sigma, \sigma \times \hat{\omega}) \in \mathbb{R}^3 \times \mathbb{R}^3 \tag{89}$$

*is a (global) constrained minimizer for $f$. Here, as in the previous section, we have denoted by $\hat{\omega}$ the unit vector associated to $\omega$.*

**Remark 18.** The construction of a couple of vectors satisfying relations (88) is straightforward. It is sufficient to choose any vector, say $\mu$, linearly independent from $\hat{\omega}$ and to set $\sigma := \mu \times \hat{\omega}$. Indeed, once defined $\hat{\mu} := \mu / |\mu|$ one has $|\sigma|^2 = |\mu|^2 (1 - (\hat{\mu} \cdot \hat{\omega})^2)$ and hence by choosing $|\mu|^2 = |\omega|/(1 - (\hat{\mu} \cdot \hat{\omega})^2)$ we get $|\sigma|^2 = |\omega|$.

**Proof.** To identify the minimizers of the problem (85)-(86), we note that for any $\omega \neq 0$, the constraint relation $v \times u = \omega$ implies the existence of a $\psi \in I_\pi := (0, \pi)$ such that $|u| \cdot |v| \cdot \sin \psi = |\omega|$. Therefore the constrained minimization for $f$ is equivalent to the unconstrained minimization of the function

$$\hat{f} : (u, \psi) \mapsto \frac{1}{2}|u|^2 + \frac{1}{2}\frac{|\omega|^2}{|u|^2 (\sin \psi)^2}, \tag{90}$$

whose stationary points $(u_\star, \psi_\star) \in \mathbb{R}^3 \times I_\pi$ are characterized by the conditions $\psi_\star = \frac{\pi}{2}$ and $|u_\star|^2 = |\omega|$. Therefore if $(u_\star, v_\star) \in \mathbb{R}^3$ is a minimizer for the function $f$ expressed by (85), then necessarily $|u_\star|^2 = |v_\star|^2 = |\omega|$, $u_\star \cdot v_\star = 0$. That this condition is also sufficient is an immediate consequence of the fact that for any such points one has $\hat{f}(u_\star, \psi_\star) = |\omega|$. Indeed, for any couple $(u, \psi) \in \mathbb{R}^3 \times I_\pi$ one has

$$\hat{f}(u, \psi) \geqslant \frac{1}{2}\frac{|u|^4 + |\omega|^2}{|u|^2} = |\omega| + \frac{1}{2}\frac{(|\omega| - |u|^2)^2}{|u|^2} \geqslant |\omega| = \hat{f}(u_\star, \psi_\star). \tag{91}$$

This concludes the proof of the statement.

For the second part we simply note that any $(u_\star, v_\star) \in \mathbb{R}^3$ can be characterized in terms of a vector orthogonal to $\omega$. Precisely, we set $\hat{\omega} := \omega/|\omega|$ and consider a vector $\sigma \in \mathbb{R}^3$ such that

$$\sigma \cdot \hat{\omega} = 0 \quad \text{and} \quad |\sigma|^2 = |\omega|. \tag{92}$$

A trivial computation shows that the vectors $v := \sigma \times \hat{\omega}$ and $u := \sigma$ satisfy the relations $u \cdot v = 0$, $|u|^2 = |v|^2 = |\omega|$ and $v \times u = \omega$. This completes the proof. $\square$

### 5.5. The characterization of the energy minimizing strokes

This section is devoted to the argument that leads to the proof of Theorem 12. It is a matter of gluing the results proved in the previous sections. We start by linking the minimization of the function $f$ to the minimization of the «$\mathfrak{g}$-orthogonalized» energy (that we rewrites here for the reader convenience, cf. (67))

$$\mathcal{G}_U(\eta) = \int_I \Lambda_{\mathfrak{g}} \dot{\eta}(t) \cdot \dot{\eta}(t) \, dt, \tag{93}$$

subject to the constraint (cf. (68)) $2\pi \langle \dot{\eta} \times \eta \rangle = \Lambda_{\mathfrak{h}}^{-1} \delta p$, where the vector $\delta p \in \mathbb{R}^3$ denotes the net displacement prescribed to the swimmer. From Proposition 17, Corollary 16 and then Proposition 14, we get that any $\sigma \in \mathbb{R}^3$ satisfying the relations

$$\sigma \cdot \hat{\omega} = 0, \quad |\sigma|^2 = |\omega|, \quad \omega := \sqrt{(\det \Lambda_{\mathfrak{g}}) \Lambda_{\mathfrak{g}}^{-1}} \Lambda_{\mathfrak{h}}^{-1} \delta p, \tag{94}$$

is associated to a (global) constrained minimizer for $\mathcal{G}_U$, via the curve $\eta(t) := (\cos t)a_1 + (\sin t)b_1$, in which the Fourier coefficients $a_1, b_1 \in \mathbb{R}^3$ are linked to the data $\omega$ (cf. (75)) by the relations $(\sqrt{2\pi\Lambda_{\mathfrak{g}}})a_1 = \sigma$ and $(\sqrt{2\pi\Lambda_{\mathfrak{g}}})b_1 = \sigma \times \hat{\omega}$. The minimum value of the energy is $\mathcal{G}_U(\eta) = |\omega|$.

Let us note that in the $\mathfrak{g}$-orthogonal reference frame the Fourier coefficients are orthonormal, i.e. once defined for every $a, b \in \mathbb{R}^3$ the inner product $(a,b)_{\mathfrak{g}} := 2\pi\Lambda_{\mathfrak{g}} a \cdot b$ and denoted by $|\cdot|_{\mathfrak{g}}$ the associated norm, the following relations hold:

$$|a_1|_{\mathfrak{g}}^2 = |b_1|_{\mathfrak{g}}^2 = |\omega| \quad \text{and} \quad (a_1, b_1)_{\mathfrak{g}} = 0. \tag{95}$$

In particular, the energy minimizing strokes for the «$\mathfrak{g}$-orthogonalized» energy $\mathcal{G}_U$, are circles in $\mathbb{R}^3$ centered at the origin and of radius $|\omega|$. Eventually, by setting $\xi := U\eta$ we get (62). This ends the proof of Theorem 12. $\square$

**Remark 19.** It is worth noticing that when the swimmer is controlled by an optimal stroke, the equations governing the dynamics of the swimmer in the range of small strokes simplify. Indeed, when $\xi$ is given by (62), from (44) we get

$$\theta(t) = \theta_0 + \nu t \quad \text{with} \quad \nu := \gamma(u \times v) \cdot \tau_3 \in \mathbb{R}, \tag{96}$$

i.e. the angular velocity of the swimmer is constant in time. In particular, when $\theta_0 = 0$ the evolution of the control system is governed by the equation (cf. (47) and (44)):

$$\dot{c}(t) = R(\nu t) F_0 \dot{\xi}(t) + R(\nu t) \sum_{j \in \mathbb{N}_2} (A_j \dot{\xi}(t) \cdot \xi(t)) \hat{e}_j. \tag{97}$$

## 6. CONCLUSIONS AND ACKNOWLEDGEMENTS

This work was partially supported by the labex LMH through the grant no. ANR-11-LABX-0056-LMH in the *Programme des Investissements d'Avenir*.